\documentclass{amsart}

\usepackage{amscd}
\usepackage[latin1]{inputenc}
\usepackage[T1]{fontenc}
\usepackage{amsmath}
\usepackage{amsfonts}
\usepackage{amssymb}
\usepackage{amsthm}
\usepackage{graphics}

\def \GL{\mathop{\rm GL}\nolimits}
\def \dim{\mathop{\rm dim}\nolimits}

\def \Fl{\mathop{\rm Fl}\nolimits}

\def \vect{\mathop{\rm vect}\nolimits}
\def \Hom{\mathop{\rm Hom}\nolimits}

\newcommand{\CC}{\mathbb{C}}
\newcommand{\KK}{\mathbb{K}}
\newcommand{\NN}{\mathbb{N}}

\newtheorem{lem}{Lemme} [section]
\newtheorem{defi}[lem]{Définition}
\newtheorem{pro}[lem]{Proposition}
\newtheorem{theo}[lem]{Théorème}
\newtheorem{cor}[lem]{Corollaire}
\newtheorem{remi}[lem]{Remarque}

\usepackage[all]{xy}
\title[espace des orbites d'une grassmannienne]{Sur la dimension de l'espace des orbites d'une grassmannienne sous l'action d'un groupe algébrique.}

\author{Michael Magen}
\address{Laboratoire de Math\'ematiques de Versailles, 
45 avenue des Etats-Unis , Bat. Fermat, 78035 VERSAILLES,
FRANCE}
\email{magen@math.math.uvsq.fr}
\urladdr{http://www.math.uvsq.fr/\textasciitilde magen}

\begin{document}

\renewcommand\contentsname{Table des matières}
\renewcommand\refname{Références}
\renewcommand\abstractname{Résumé}
%\pagestyle{myheadings}
%\markboth{MICHAEL MAGEN}{ESPACE DES ORBITES D'UNE GRASSMANNIENNE}

\begin{abstract}

\vskip 4.5mm

On s'intéresse aux actions d'un groupe algébrique $G$ sur les grassmanniennes d'un $\KK$-espace vectoriel de dimension finie $V$ ($\KK$ un corps algébriquement clos) déduites d'une action de $G$ sur $V$. On montre que \linebreak $\dim_G G(j,V) \leq \dim_G G(k,V)$ si et seulement si  $\dim G(j,V) \leq \dim G(k,V)$ où $\dim_G G(r,V)$ désigne la dimension de l'espace des orbites de $G(r,V)$ sous l'action de $G$ (ce qui généralise un résultat de Pyasetskii {\cite{PYA}}). On étend ensuite ce résultat aux variétés de drapeaux. Des méthodes différentes nous permettent d'obtenir des résultats sur les nombres d'orbites , quand le corps de base $\KK$ est fini.

\noindent {\bf 2000 Mathematics Subject Classification:}  20G15, 20G40, 14L30.

\end{abstract}

\maketitle

\tableofcontents

\section{Introduction}
L'objet de cet article est de généraliser le résultat suivant dû à Pyasetskii {\cite{PYA}} :\\

\begin{theo}
~\\
Soit $G$ un groupe algébrique agissant sur un $\CC$-espace vectoriel $V$ de dimension finie $n$. Soit $\mathbb{P}(V)$ et $\mathbb{P}(V^*)$ les espaces projectifs paramétrant respectivement les droites et les hyperplans de $V$, tous deux naturellement munis d'une action de $G$. Les deux assertions suivantes sont équivalentes :
$$\begin{array}{c}
(i)\text{ Le groupe $G$ agit avec un nombre fini d'orbites sur $\mathbb{P}(V)$},\\
(ii)\text{ Le groupe $G$ agit avec un nombre fini d'orbites sur $\mathbb{P}(V^*)$}.\\
\end{array}
$$ 
\end{theo}

Désormais, sauf dans les sections \ref{fini} et \ref{fini2}, on suppose le corps de base $\KK$ algébriquement clos.\\

L'assertion de Pyasetskii suggère une assertion analogue pour les grassmanniennes de taille $k$ et $n-k$ mais on n'en trouve aucune démonstration  (à notre connaissance) dans la littérature. Nous allons montrer en particulier qu'une telle assertion est vraie (corollaire \ref{Pya2}).\\

\begin{defi}
~\\
Soit $G$ un groupe algébrique agissant sur une variété algébrique $X$. D'après le théorème de Rosenlicht \cite{ROS}, on peut trouver une partition de $X$ : $X=\coprod\limits_{i=1}^{r}X_i$ tel que les $X_i$ sont des sous-ensembles constructibles $G$-invariants de $X$ admettant tous un quotient géométrique. On note $\dim_G X$ le maximum des dimensions de ces quotients (ce nombre est bien défini, i.e. indépendant de la partition ad hoc utilisée).\\
\end{defi}

Pour $V$ un $ \KK$-espace vectoriel de dimension $n$ et pour $k$ dans $\{ 1,\dots,n-1 \}$, on note $G(k,V)$ la grassmannienne des sous-espaces de dimension $k$ de $V$.\\

Soit $G$ un groupe algébrique agissant sur un $\KK$-espace vectoriel $V$ de dimension finie $n$. Pour tout $j$ dans $\{ 1,\dots,n-1 \}$, la grassmannienne $G(j,V)$ est naturellement munie d'une action de $G$. Le résultat principal de l'article est le suivant :

\begin{theo} \label{main}
~\\
On a :
$$\begin{array}{c}
\dim_G G(j,V) \leq \dim_G G(k,V) \text{ pour tout couple $(j,k)$ tel que $\dim G(j,V) \leq \dim G(k,V)$}. \\
\end{array}
$$
\end{theo}

On a en particulier les deux résultats suivants :\\

\begin{cor}
~\\
Pour tout $k$ dans $\{ 1,\dots,n-1 \}$, on a : $$\dim_G G(k,V)=\dim_G G(n-k,V).$$
\end{cor}
~

\begin{cor} \label{Pya2}
~\\
Pour tout $k$ dans $\{ 1,\dots,n-1 \}$, les assertions suivantes sont équivalentes :
$$\begin{array}{cl}
(i) & \text{Le groupe $G$ agit avec un nombre fini d'orbites sur $G(k,V)$} \\
(ii) & \text{Le groupe $G$ agit avec un nombre fini d'orbites sur $G(j,V)$}\\
     & \text{pour tout $j$ tel que $\dim G(j,V) \leq \dim G(k,V)$}. \\
\end{array}
$$
\end{cor}

Pour démontrer le théorème \ref{main}, on utilise une généralisation géométrique de la formule de Burnside (théorème \ref{burnside}). On démontre un résultat analogue pour les variétés de drapeaux dans la section \ref{secdrap} et dans la section \ref{secfini}, ces résultats sont adaptés au contexte d'un corps de base fini.

\section{Grassmanniennes, variétés de drapeaux et variétés de Burnside }

Pour toutes les notions concernant les partitions et les compositions d'un entier utilisées dans cet article (conjugaison, "raising operators",...), l'ouvrage de référence est celui de Macdonald \cite{McD}.\\

\subsection{Variétés de drapeaux} 
~\\
Soit $n$ un entier strictement positif.\\
\begin{defi}
~\\
Soit $k$ un entier strictement positif. On dit que $\underline{a}=(a_1,\dots ,a_k)\in  (\NN)^k$ est une composition de $n$ de longueur $k$ si $a_1+\dots +a_k=n$.\\
Les nombres $a_1,\dots ,a_k$ sont appelés les termes de $\underline{a}$.
\end{defi}

Soit $V$ un $\KK$-espace vectoriel de dimension finie $n$.\\

\begin{defi}
~\\
Pour toute composition $(a_1,\dots ,a_r)$ de n, on définit la variété de drapeaux :

$$\Fl_{(a_1,\dots ,a_r)}(V)=
\left\lbrace
\begin{array}{c}
0=V_0 \subset V_1 \subset \dots V_r=V ,\\
\text{les $V_i$ sont des sous-espaces vectoriels de $V$ et} \\
\text{pour tout $i$, } \dim V_i/V_{i-1}=a_i \\
\end{array}
\right\rbrace
.
$$
\end{defi}

Si $G$ est un groupe algébrique agissant sur $V$, $\Fl_{(a_1,\dots ,a_r)}(V)$ hérite naturellement d'une action de $G$.\\

\begin{remi}
~\\
Pour tout $k$ dans $\{ 1,\dots,n-1 \}$, on a $G(k,V)=\Fl_{(k,n-k)}(V)$.
\end{remi}
~

\subsection{La variété de Burnside d'une $G$-variété algébrique}
~\\

Soit $G$ un groupe algébrique et $X$ une $G$-variété.\\

On définit $\Sigma$, la variété de Burnside de $X$ sous l'action de $G$, de la manière suivante :
$$ \Sigma= \{ (g,x)\in G \times X \ / \ g.x=x \}.$$
On a $\dim \Sigma \geq \dim G$.

\begin{theo}\label{burnside}
~\\
On a l'égalité suivante :
$$\dim_G X= \dim \Sigma -\dim G.$$
\end{theo}

{\sl Démonstration~:}\\
On peut supposer que le quotient $X/G$ est géométrique et qu'il existe un entier positif $k$ tel que pour tout $x$ dans $X$, $\dim G.x=k$.
On considère le morphisme de projection de $\Sigma$ dans $X$, la fibre d'un élément de $x$ est isomorphe à $\{g\in G \ , \ g.x=x \}$ et est donc de dimension constante $\dim G-k$. On a donc $\dim \Sigma =\dim X + \dim G -k$.\\
Par ailleurs, d'après le théorème de Rosenlicht, on a $\dim (X/G)=\dim X -k$ d'où le résultat.$\Box$\\

\begin{remi}
~\\
Le théorème \ref{burnside} fournit une nouvelle démonstration du fait que le nombre $\dim_G X$ est bien défini.
\end{remi}

Soit $g$ dans $G$, on note $X_g$ l'ensemble des points de $X$ fixés par $g$.\\

\begin{remi}
~\\
Soit $q:\Sigma \rightarrow G$ le morphisme de projection. Pour tout $g$ dans $G$, la fibre $q^{-1}(\{g\})$ est isomorphe à $X_g$.
\end{remi}

D'après le théorème \ref{burnside}, le calcul de la dimension de l'espace des orbites de l'action de $G$ sur $X$ se ramène au calcul de la dimension de la variété de Burnside. Le calcul de cette dernière dimension peut se faire  en étudiant le morphisme de projection de  $\Sigma$ dans $G$, en utilisant une partition finie de $G$ en constructibles sur lesquels la dimension de la fibre est constante.\\

Dans la section suivante, on construit une telle partition pour $\GL(n,\KK)$ quand $X$ est une variété de drapeaux.\\

\subsection{Partitions d'un entier, classes de conjugaison unipotentes et squelettes dans $\GL$}
~\\

Soit $n$ un entier strictement positif.\\

\begin{defi}
~\\
Soit $k$ un entier strictement positif. On dit que $\lambda=(\lambda_1,\dots ,\lambda_k)\in  (\NN)^k$ est une partition de $n$ de longueur k si $\lambda_1+\dots +\lambda_k=n$ et $\lambda_1 \geq \dots \geq \lambda_k \geq 0$.\\
Les nombres $\lambda_1,\dots ,\lambda_k$ sont appelés les termes de $\lambda$.
\end{defi}

On note $\mathcal{P}(n)$ l'ensemble des partitions de $n$. Rappelons la définition de l'ordre de dominance sur $\mathcal{P}(n)$.\\

\begin{defi}
~\\
Soit $\lambda=(\lambda_1,\dots ,\lambda_r)$ et $\mu=(\mu_1,\dots ,\mu_r)$ dans $\mathcal{P}(n)$. On note $\lambda \geq \mu$ si pour tout $i$ dans $\{ 1,\dots,r \}$, $\lambda_1 +\dots + \lambda_i \geq \mu_1 + \dots + \mu_i$. 
\end{defi}
~

\begin{defi}
~\\
Soit $m$ un entier strictement positif et $z$ un élément de $\KK$. On pose :
$$ J_m(z)=
\left(
\begin{array}{ccccc}
z      & 1 & 0 & \dots & 0\\
0      & z & 1 & \dots & 0\\
\vdots &   &   &       & \vdots \\
0      & 0 & 0 & \dots & z\\
\end{array}
\right).
 $$
Soit $\lambda=(\lambda_1,\dots ,\lambda_r)$ une partition de $n$ dont tous les termes sont strictement positifs et $z$ un élément de $\KK$. On pose :
$$ J_{\lambda}(z)=
\left(
\begin{array}{ccccc}
J_{\lambda_1}(z)      & 0 & 0 & \dots & 0\\
0      & J_{\lambda_2}(z) & 0 & \dots & 0\\
\vdots &   &   &       & \vdots \\
0      & 0 & 0 & \dots & J_{\lambda_r}(z)\\
\end{array}
\right).
 $$
\end{defi}

On rappelle la proposition suivante :\\
\begin{pro}
~\\
Soit $u$ un unipotent dans $\GL(n,\KK)$, il existe une unique partition $\lambda$ de $n$ dont tous les termes sont strictement positifs telle que $u$ soit conjugué à $J_{\lambda}(1)$. On dit alors que $u$ est un unipotent de type $\lambda$.
\end{pro}
~

\begin{defi}
~\\
Soit $u$ et $v$ dans $\GL(n,\KK)$, on note $ \text{ $u \ \mathcal{R} \ v$}$ s'il existe un entier $r$, deux $r$-uplets $(a_1,\dots ,a_r)$ et $(b_1,\dots ,b_r)$ d'éléments de $\KK$ deux à deux distincts et $r$ partitions $\lambda^{(1)},\dots ,\lambda^{(r)}$ tels que :
$$u \text{ est conjugué à }
\left(
\begin{array}{ccccc}
J_{\lambda^{(1)}}(a_1)      & 0 & 0 & \dots & 0\\
0      & J_{\lambda^{(2)}}(a_2) & 0 & \dots & 0\\
\vdots &   &   &       & \vdots \\
0      & 0 & 0 & \dots & J_{\lambda^{(r)}}(a_r)\\
\end{array}
\right),
$$
$$
\text{ et $v$ est conjugué à }
\left(
\begin{array}{ccccc}
J_{\lambda^{(1)}}(b_1)      & 0 & 0 & \dots & 0\\
0      & J_{\lambda^{(2)}}(b_2) & 0 & \dots & 0\\
\vdots &   &   &       & \vdots \\
0      & 0 & 0 & \dots & J_{\lambda^{(r)}}(b_r)\\
\end{array}
\right).
$$
\end{defi}

Clairement, la relation $\mathcal{R}$ est une relation d'équivalence sur $\GL(n,\KK)$. On note $\mathcal{S}(n)$ l'ensemble de ses classes d'équivalence ($\mathcal{S}(n)$ est fini). Pour $s$ dans $\mathcal{S}(n)$, on note $\GL^{(s)}(n,\KK)$ l'ensemble des éléments de $\GL(n,\KK)$ dans la classe d'équivalence $s$.\\

\begin{defi}
~\\
Les éléments de $\mathcal{S}(n)$ sont appelés squelettes de taille $n$.
\end{defi}
~
\begin{defi}
~\\
Un élément $s$ de $\mathcal{S}(n)$ est dit unipotent (respectivement semisimple) s'il existe un élément unipotent (respectivement semisimple) dans $\GL^{(s)}(n,\KK)$.
\end{defi}
~

\begin{pro}
~\\
Pour tout $s$ dans $\mathcal{S}(n)$, $\GL^{(s)}(n,\KK)$ est constructible. Ainsi, on a la partition suivante de $\GL(n,\KK)$ en constructibles :
$$ \GL(n,\KK)=\coprod\limits_{s \in \mathcal{S}(n)} \GL^{(s)}(n,\KK).$$
\end{pro}

Nous donnons à titre d'exemple cette partition pour $n=3$ : 

$$\begin{array}{llll}
\GL(3,\KK) & = & & 
\left\lbrace
\left(
\begin{array}{ccc}
a & 0 & 0 \\
0 & a & 0 \\
0 & 0 & a \\
\end{array}
\right)
, \ a \in \KK^*.
\right\rbrace \\
\\
 & & \cup &
\left\lbrace
h
\left(
\begin{array}{ccc}
a & 1 & 0 \\
0 & a & 0 \\
0 & 0 & a \\
\end{array}
\right)
h^{-1}
,
\left\lbrace 
\begin{array}{l}
a \in \KK^*, \\
h \in \GL(3,\KK).\\
\end{array}
\right.
\right\rbrace \\
\\
 & & \cup & 
\left\lbrace
h
\left(
\begin{array}{ccc}
a & 1 & 0 \\
0 & a & 1 \\
0 & 0 & a \\
\end{array}
\right)
h^{-1}
,
\left\lbrace 
\begin{array}{l}
 a \in \KK^* ,\\
 h \in \GL(3,\KK).\\
\end{array}
\right.
\right\rbrace \\
\\
& & \cup &
\left\lbrace
h
\left(
\begin{array}{ccc}
a & 0 & 0 \\
0 & a & 0 \\
0 & 0 & b \\
\end{array}
\right)
h^{-1}
,
\left\lbrace 
\begin{array}{l}
 (a,b) \in (\KK^*)^2 ,\\
a \neq b ,\\
h \in \GL(3,\KK).\\
\end{array}
\right.
\right\rbrace \\
\\
 & & \cup & 
\left\lbrace
h
\left(
\begin{array}{ccc}
a & 1 & 0 \\
0 & a & 0 \\
0 & 0 & b \\
\end{array}
\right)
h^{-1}
,
\left\lbrace 
\begin{array}{l}
(a,b) \in (\KK^*)^2 ,\\
a \neq b ,\\
h \in \GL(3,\KK).\\
\end{array}
\right.
\right\rbrace \\
\\
& & \cup &
\left\lbrace
h
\left(
\begin{array}{ccc}
a & 0 & 0 \\
0 & b & 0 \\
0 & 0 & c \\
\end{array}
\right)
h^{-1}
,
\left\lbrace 
\begin{array}{l}
(a,b,c) \in (\KK^*)^3 ,\\
a \neq b \neq c \neq a,\\
h \in \GL(3,\KK).\\
\end{array}
\right.
\right\rbrace \\
\end{array}
$$

Nous allons maintenant énoncer certaines propriétés remarquables de cette partition, propriétés dont nous allons nous servir de manière essentielle dans la démonstration du résultat principal.

\begin{pro} \label{fibre}
~\\
Soit $\underline{a}$ une composition de  $n$. La fonction :
$$\left\lbrace
\begin{array}{lll}
\GL(n,\KK) & \rightarrow & \NN \\
g      & \mapsto     & \dim (\Fl_{\underline{a}}(V))_g \\
\end{array}
\right.
$$
est constante sur les $\GL^{(s)}(n,\KK)$.
\end{pro}

\begin{pro} \label{pro2}
~\\
Soit $g$ un élément de $\GL(n,\KK)$, il existe un squelette $t$ semisimple tel que pour tout $v$ dans $\GL^{(t)}(n,\KK)$ et pour toute composition $\underline{a}$ de $n$ : 
$$\dim (\Fl_{\underline{a}}(V))_g =\dim (\Fl{\underline{a}}(V))_v.$$
\end{pro}

\begin{pro} \label{dual}
~\\
Soit $g$ un élément de $\GL(n,\KK)$, on a :
$$\dim G(j,V)_g  \leq G(k,V)_g \ \text{pour tout couple $(j,k)$ tel que $\dim G(j,V) \leq \dim G(k,V)$}.$$
\end{pro}

Cette dernière assertion donne les deux résultats suivants dans le cadre des va\-rié\-tés de drapeaux :

\begin{cor} \label{dualdrap}
~\\
Soit $g$ un élément de $\GL(n,\KK)$, et $(a_1,\dots ,a_r)$ une composition de $n$. On a :
$$\dim (\Fl_{(a_1,\dots ,a_r)}(V))_g=\dim (\Fl_{(a_{\sigma(1)},\dots ,a_{\sigma(r)})}(V))_g \ \text{pour tout $\sigma \in \mathfrak{S}_r$}, $$
où $\mathfrak{S}_r$ désigne le groupe des permutations de l'ensemble $\{1,\dots,r\}$.
\end{cor}

\begin{cor} \label{dualdrap2}
~\\
Soit $g$ un élément de $\GL(n,\KK)$, $\underline{a}=(a_1,\dots ,a_r)$ une composition de $n$, $i$ dans $\{ 1,\dots,r-1 \}$ et $b$ et $c$ deux entiers tels que $b+c=a_i+a_{i+1}$ et $bc \leq a_ia_{i+1}$ (i.e $\dim G(b,b+c) \leq \dim G(a_i,a_i+a_{i+1})$). On a :
$$ \dim (\Fl_{(a_1,\dots ,a_{i-1},b,c,a_{i+2},\dots a_r)}(V))_g \leq \dim (\Fl_{(a_1,\dots ,a_r)}(V))_g.$$
\end{cor}

La proposition \ref{fibre} est classique, démontrons les quatre autres assertions.\\

{\sl Démonstration de la proposition \ref{pro2}~:}\\
On se ramène facilement (par la considération des sous-espaces caractéristiques) au cas où $g$ est unipotent.\\
Supposons donc que $x$ est unipotent de type $\lambda$. Soit $y$ un élément semisimple de $\GL(n,\KK)$ et $\mu=(\mu_1,\dots ,\mu_s)$ une partition de $n$. On dit que $y$ est de type $\mu$ s'il existe $s$ éléments non nuls de $\KK$ deux à deux distincts $z_1,\dots ,z_s$ tels que :
$$y=\left(
\begin{array}{ccc}
z_1Id_{\mu_1} & & \\
 & \ddots & \\
 & & z_sId_{\mu_s} \\
\end{array}
\right).
$$ 
On note $\lambda'$ la partition conjuguée de $\lambda$, et soit $t$ le squelette semisimple tel que $\GL^{(t)}(n,\KK)$ est formé des semisimples de type $\lambda'$. Soit $v$ dans $\GL^{(t)}(n,\KK)$.\\

On pose $b_i=a_1+\dots +a_i$ pour tout $i$ dans $\{ 1,\dots,r \}$, $\mathfrak{M}=\Hom(\KK^{b_1},\KK^{b_2}) \oplus \dots \oplus \Hom(\KK^{b_{r-1}},\KK^{b_r})$ et $G=\GL(b_1,\KK) \times \dots \times \GL(b_r,\KK)$.On se donne une suite d'injections $\KK^{b_1} \hookrightarrow \KK^{b_2} \hookrightarrow \dots \hookrightarrow \KK^{b_r}$ de sorte que l'on a une flèche de restriction $\GL(b_r,\KK) \rightarrow G$. On définit maintenant les deux applications suivantes :
$$ \phi :\left\{
\begin{array}{ccc}
(\Fl_{\underline{a}})_x & \rightarrow & \mathcal{P}(b_1) \times \dots \times \mathcal{P}(b_r) \\
(V_0 \subset \dots \subset V_r) & \mapsto & (\lambda^{(1)},\dots ,\lambda^{(r)})\\
\end{array}
\right., 
$$
où $\lambda^{(i)}$ est le type de la restriction de $x$ à $V_i$ pour tout $i$ dans $\{ 1,\dots,r \}$ et

$$ \psi :\left\{
\begin{array}{ccc}
(\Fl_{\underline{a}})_v & \rightarrow & \mathcal{P}(b_1) \times \dots \times \mathcal{P}(b_r) \\
(V_0 \subset \dots \subset V_r) & \mapsto & (\lambda^{(1)},\dots ,\lambda^{(r)})\\
\end{array}
\right.,
$$
où $\lambda^{(i)}$ est le type de la restriction de $v$ à $V_i$ pour tout $i$ dans $\{ 1,\dots,r \}$.

On obtient donc les deux partitions (en ensembles constructibles) suivantes des ensembles $(\Fl_{\underline{a}})_x$ et $(\Fl_{\underline{a}})_v$ :
$$\begin{array}{c}
(\Fl_{\underline{a}})_x=\coprod\limits_{\alpha \in \mathcal{P}(b_1) \times \dots \times \mathcal{P}(b_r)} \phi^{-1}(\{ \alpha\}), \\
(\Fl_{\underline{a}})_v=\coprod\limits_{\beta \in \mathcal{P}(b_1) \times \dots \times \mathcal{P}(b_r)} \psi^{-1}(\{ \beta\}). \\ 
\end{array}
$$

Nous allons montrer que pour tout $\alpha$ dans $\mathcal{P}(b_1) \times \dots \times \mathcal{P}(b_r)$, $\dim \phi^{-1}(\{ \alpha\})= \dim \psi^{-1}(\{ \alpha'\})$ , où $\alpha'$ est obtenue en prenant les conjuguées des termes de $\alpha$.\\

Soit donc $\alpha=(\alpha^{(1)},\dots ,\alpha^{(r)})$ dans $\mathcal{P}(b_1) \times \dots \times \mathcal{P}(b_r)$. On pose $\underline{x}=(x_1,\dots ,x_r)$ dans $G$ où $x_i$ est la restriction de $x$ à $\KK^{b_i}$ ($x_i$ est un unipotent de type $\alpha^{(i)}$) pour tout $i$ dans $\{ 1,\dots,r \}$. On pose aussi $\underline{v}=(v_1,\dots ,v_r)$ dans $G$ où $v_r$ est conjugué à $v$ et $v_i$ est la restiction de $v$ à $\KK^{b_i}$ et est un semisimple de type la conjuguée de $\alpha^{(i)}$ pour tout $i$ dans $\{ 1,\dots,r \}$.\\

D'après \cite[lemme 3]{KAC}, on a :
$$ \dim \mathfrak{M}^{\underline{x}}=\dim \mathfrak{M}^{\underline{v}}.$$

Par ailleurs, on a deux flèches surjectives :
$$f:\left\lbrace
\begin{array}{ccc}
\mathfrak{M}^{\underline{x}} & \rightarrow & \phi^{-1}(\{ \alpha\})\\
(g_1,\dots ,g_{r-1}) & \mapsto & (g_1(\KK^{b_1}) \subset \dots \subset g_{r-1}(\KK^{b_{r-1}}))\\
\end{array}
\right.,
$$ 
et 
$$g:\left\lbrace
\begin{array}{ccc}
\mathfrak{M}^{\underline{v}} & \rightarrow & \psi^{-1}(\{ \alpha'\})\\
(g_1,\dots ,g_{r-1}) & \mapsto & (g_1(\KK^{b_1}) \subset \dots \subset g_{r-1}(\KK^{b_{r-1}}))\\
\end{array}
\right..
$$ 

Pour tout $\underline{V}$ dans $\phi^{-1}(\{ \alpha\})$, $f^{-1}(\{ \underline{V} \})$ est isomorphe à $\GL(b_1,\KK)^{x_1} \times \dots \times \GL(b_r-1,\KK)^{x_{r-1}}$ et pour tout $\underline{W}$ dans $\psi^{-1}(\{ \alpha'\})$, $g^{-1}(\{ \underline{W} \})$ est isomorphe à $\GL(b_1,\KK)^{v_1} \times \dots \times \GL(b_r-1,\KK)^{v_{r-1}}$, or ces deux espaces ont même dimension d'où le résultat.$\Box$\\

{\sl Démonstration de la proposition \ref{dual}~:}\\
La première assertion est un simple résultat de dualité. Pour démontrer la seconde, on observe que d'après la proposition \ref{pro2} on peut supposer que $x$ est semisimple. Le résultat est alors élémentaire.$\Box$\\

{\sl Démonstration du corollaire \ref{dualdrap}~:}\\
Soit $i$ dans $\{ 1,\dots,r-1 \}$, montrons que le résultat est vrai pour pour la transposition $\tau=(i,i+1)$. Notons $X=\Fl_{(a_1,\dots ,a_r)}(V)$, $Y=\Fl_{(a_{\tau(1)},\dots ,a_{\tau(r)})}(V)$ et $Z=\Fl_{(a_1,\dots,a_{i-1},a_i+a_{i+1},a_{i+2},\dots ,a_r)}(V)$. On définit les flèches suivantes :
$$\phi :\left\{
\begin{array}{ccc}
X_g & \rightarrow & Z_g\\
(V_0 \subset \dots \subset V_r) & \mapsto & (V_0 \subset \dots \subset V_{i-1} \subset V_{i+1} \subset V_r)\\
\end{array}
\right. 
$$
et
$$\psi :\left\{
\begin{array}{ccc}
Y_g & \rightarrow & Z_g\\
(V_0 \subset \dots \subset V_r) & \mapsto & (V_0 \subset \dots \subset V_{i-1} \subset V_{i+1} \subset V_r)\\
\end{array}
\right..
$$

Soit $\underline{V}=(V_0 \subset \dots \subset V_{i-1} \subset V_{i+1} \subset V_r)$ dans $Z_g$, on note $\tilde{g}$ l'élément de $\GL(V_{i+1}/V_{i-1})$ déduit de $g$. On a $u^{-1}(\{\underline{V}\})=G(a_i,V_{i+1})_{\tilde{g}}$ et $v^{-1}(\{\underline{V}\})=G(a_{i+1},V_{i+1})_{\tilde{g}}$. Donc d'après la proposition \ref{dual}, pour tout $ \underline{V}$ dans $Z_g$, $\dim u^{-1}(\{\underline{V}\})=\dim v^{-1}(\{\underline{V}\})$, et $\phi$ et $\psi$ sont surjectives, d'où le résultat.$\Box$\\

Le corollaire \ref{dualdrap2} se montre de la même manière. \\

\subsection{Démonstration du théorème \ref{main}}
~\\

Le groupe $G$ agissant sur $V$, on peut supposer que c'est un sous-groupe de $\GL(V)=\GL(n,\KK)$. Pour tout $j$ dans $\{ 1,\dots,n \}$, on note $\Sigma_j$ la variété de Burnside de $G(j,V)$ :
$$\begin{array}{c}
\Sigma_j = \{ (g,P) \in G \times G(j,V) \ / \ g.P=P \}.\\
\end{array}
$$
On note $p_j$ la projection de $\Sigma_j$ dans $G$. On a : $ G=\coprod\limits_{s \in \mathcal{S}(n)} (\GL^{(s)}(n,\KK) \cap G)$. Et donc $\Sigma_j =\coprod\limits_{s \in \mathcal{S}(n)} p_j^{-1}(\GL^{(s)}(n,\KK) \cap G)$. Puis $\dim \Sigma_j=\text{max } \{ \dim p_j^{-1}(\GL^{(s)}(n,\KK) \cap G) \}$. Pour $g \in G$, $p_j^{-1}(\{ g \})=\{P \in G(j,V) \ / \ g.P=P \}$, par conséquent d'après la proposition \ref{fibre}, la dimension de $p_j^{-1}(\{g\})$ est constante sur les $\GL^{(s)}(n,\KK) \cap G$ pour tout $s$ dans $\mathcal{S}(n)$, notons la $d_j(s)$. On a donc $\dim p_j^{-1}(\GL^{(s)}(n,\KK) \cap G)=\dim (\GL^{(s)}(n,\KK) \cap G)+d_j(s)$.\\
 
D'après la proposition \ref{dual}, $d_j(s) \leq d_{k}(s)$ pour tout $s \in \mathcal{S}(n)$ et tout $j \in E_k$. On obtient donc $\dim \Sigma_j \leq \dim \Sigma_{k}$ ce qui nous permet de conclure d'après le théorème \ref{burnside}.$\Box$

\subsection{Un résultat analogue pour les variétés de drapeaux} \label{secdrap}
~\\

Soit $G$ un groupe algébrique et $V$ un $\KK$-espace vectoriel de dimension finie $n$ muni d'une action de $G$.\\

Soit $(a_1,\dots ,a_r)$ une composition de $n$. Il existe une permutation $\sigma$ dans $\mathfrak{S}_r$ tel que $(a_{\sigma(1)},\dots ,a_{\sigma(r)})$ est une partition de $n$. On note $P(a_1,\dots ,a_r)$ cette partition.\\

On démontre comme on l'a fait pour les grassmanniennes les résultats suivants :

\begin{theo} \label{drapeaux}
~\\
Soit $(a_1,\dots ,a_r)$ une composition de $n$. On a l'égalité suivante :
$$\dim_G \Fl_{(a_1,\dots ,a_r)}(V)=\dim_G \Fl_{(a_{\sigma(1)},\dots ,a_{\sigma(r)})}(V) \ \text{pour tout $\sigma \in \mathfrak{S}_r$}. $$
\end{theo}

\begin{pro} \label{bleu} 
~\\
Soit $\underline{a}=(a_1,\dots ,a_r)$ une composition de $n$, $i$ dans $\{ 1,\dots,r-1 \}$ et $b$, $c$ deux entiers tels que $b+c=a_i+a_{i+1}$ et $bc \leq a_ia_{i+1}$, i.e $\dim G(b,b+c) \leq \dim G(a_i,a_i+a_{i+1})$. On a :
$$ \dim_G \Fl_{(a_1,\dots ,a_{i-1},b,c,a_{i+2},\dots a_r)}(V) \leq \dim_G \Fl_{(a_1,\dots ,a_r)}(V).$$
\end{pro}

On déduit de ces deux assertions le résultat suivant :\\

\begin{theo}
~\\
Soit $\underline{a}=(a_1,\dots ,a_r)$ et $\underline{b}=(b_1,\dots ,b_r)$ deux compositions de $n$ telles que $P(\underline{a}) \geq P(\underline{b})$. On a : 
$$ \dim_G \Fl_{\underline{a}}(V) \leq \dim_G \Fl_{\underline{b}}(V).$$
\end{theo}

Ce théorème se déduit directement des deux résultats précédents et du lemme \ref{comb} (ce lemme est démontré en annexe et est totalement indépendant du reste de l'article).

\section{Grassmaniennes et variétés de drapeaux sur un corps fini} \label{secfini}
\subsection{Grassmaniennes sur un corps fini} \label{fini}
~\\

Soit $\mathbb{F}_q$ un corps fini à $q$ éléments. Soit maintenant $n$ un entier naturel non nul et $V$ un $\mathbb{F}_q$-espace vectoriel de dimension finie $n$.\\

Soit $E$ un ensemble fini, on note $\# E$ son cardinal. Si $E$ est muni d'une action d'un groupe $G$, on note $N(G,E)$ le nombre d'orbites de cette action.

\begin{pro} \label{grassfin}
~\\
Soit $G$ un groupe agissant sur $V$. On déduit de cette action une action sur $G(r,V)$ pour tout $r$ dans $\{ 1,\dots,n-1 \} $. On a : 
$$  N(G,G(j,V)) \leq N(G,G(k,V)) \text{pour tout couple $(j,k)$ tel que $\# G(j,V) \leq \# G(k,V)$.}$$
\end{pro}

{\sl Démonstration~:}\\
On commence par remarquer que le résultat annoncé est équivalent au résultat suivant : pour tout $k$ dans $\{ 1,\dots,n/2 \}$, on a les deux assertions suivantes :
$$\begin{array}{ll}
(i) & N(G,G(k,V)) \geq N(G,G(r,V)) \text{ pour tout $r$ dans $\{ 1,\dots,k \}$},\\
(ii) & N(G,G(k,V))=N(G,G(n-k,V)) \\
\end{array}
$$

Démontrons l'assertion $(i)$.\\

Pour $r$ dans $\{ 1,\dots,n-1 \}$, on note $F_r$ le $\KK$-espace vectoriel de dimension $\#G(r,V)$ des fonctions de $G(r,V)$ dans $\KK$. Le sous-espace vectoriel de $F_r$ formé par les fonctions $G$-invariantes est noté $F_r^G$. On a l'égalité :
$$ N(G,G(r,V))=\dim F_r^G. $$ \\

Soit donc $k$ dans $\{ 1,\dots,n/2 \}$, et $r$ dans $\{ 1,\dots,k \}$

On définit l'application suivante :\\
$$ \phi :\left\{
\begin{array}{ccc}
F_r & \rightarrow & F_k \\
f & \mapsto & \hat{f}\\
\end{array}
\right.,$$ 
où $\hat{f}$ est définie par :\\
$$ \hat{f} :\left\{
\begin{array}{ccc}
G(k,V) & \rightarrow & \mathbb{C} \\
H & \mapsto & \displaystyle{\sum_{P \subset H} f(P)} \\
\end{array}
\right..$$ 

Cette application est $G$-équivariante, nous allons montrer qu'elle est injective, la première partie du résultat sera alors démontré.\\

Pour cela, on va construire une application $\pi : F_k \rightarrow F_r $ telle que $\pi \circ \phi=Id_{F_r}$.\\
On pose :

$$ \pi :\left\{
\begin{array}{ccc}
F_{k} & \rightarrow & F_{r} \\
g & \mapsto & \check{g}\\
\end{array}
\right.,$$ 
où $\check{g}$ est définie par :\\
$$ \check{g} :\left\{
\begin{array}{ccc}
G(r,V) & \rightarrow & \mathbb{C} \\
P & \mapsto & \displaystyle{\sum_{H}} \epsilon_P(H)g(H) \\
\end{array}
\right.,$$ 
et les $\epsilon_P(H)$ sont des nombres complexes. L'exercice consiste à montrer que l'on peut définir ces coefficients de manière à obtenir l'égalité $\pi \circ \phi=Id_{F_k}$.\\

Soit $P_0\in G(r,V)$ , on d\'efinit  $f_{P_0}\in F_{r}$ par :\\
$$f_{P_0}(P)=\left\{
\begin{array} {cc}
1 & \text{si } P =P_0 \\
0 & \text{sinon}.\\
\end{array}
\right.$$
 
On a donc :\\
$$\forall H \in G(k,V) \quad \hat{f}_{P_0}(H)=\left\{
\begin{array} {cc}
1 & \text{si } P \subset H \\
0 & \text{sinon}. \\
\end{array}
\right.$$

Puis :\\
$$ \begin{array}{cl} 
\forall P \in G(r,V) \quad \check{\hat{f}}_{P_0}(P) &  
=\displaystyle{\sum_{H}} \epsilon_P(H)\hat{f}_{P_0}(H)\\
& =\displaystyle{\sum_{H\supset P_0}} \epsilon_P(H).\\
\end{array}
$$

On veut donc :\\
$$ \displaystyle{\sum_{H\supset P_0}} \epsilon_P(H)=
\left\{
\begin{array}{cc}
0 & si \quad P\neq P_0,\\
1 & si \quad P=P_0.\\
\end{array}
\right.$$

Pour ce faire, on va supposer que $\epsilon_P(H)$ ne d\'epend que de la dimension de l'intersection de $P$ et $H$. On pose donc $\epsilon_P(H)=\epsilon_{\dim P\cap H}$. Il nous faut donc d\'efinir les $\epsilon_j$ pour $ j \in \{ 0,\dots,r \} $ de manière convenable.\\

Soit $i \in \{ 0,\dots,r \}$  et soit $P$ et $P_0$ dans $G(r,V)$ tels que $\dim(P \cap P_0)=i$. On pose pour $j \in \{ 0,\dots,r \}$, $a_{i,j}=\# \{ H \in G(k,V)/P_0 \subset H \, , \, \dim(P \cap H)=j \}$. Ces nombres ne d\'ependent pas des choix de $P$ et $P_0$. Soit $A$ la matrice dont les coefficients sont les $a_{i,j}$. Montrons que $A$ est inversible. On a $a_{i,j}=0$ pour $i>j$ (car $P \cap P_0 \subset P \cap H $ ) i.e. A est triangulaire sup\'erieure.\\

Montrons maintenant que pour tout $i \in \{ 0,\dots,r \}$, $a_{i,i}$ est non nul. Soit donc $i \in \{ 0,\dots,r \}$, on doit trouver un $H$ dans $G(k,V)$ tel que $P_0 \subset H$ et $\dim P \cap H=i$. Nous allons en construire un.
On choisit une base $(e_1,\dots ,e_n)$ de $V$ telle que $P_0=\vect(e_1,\dots ,e_r)$ et $P_1=\vect(e_1,\dots ,e_i,e_{r+1},\dots ,e_{2r-i})$. La famille $(e_{2r-i+1},\dots ,e_n)$ comporte $n-2r+i$ éléments et on a l'inégalité suivante : $k-r \leq n-2r+i$. Par conséquent, on peut poser $$H=\vect(e_1,\dots ,e_r,e_{2r-i+1},\dots ,e_{k+r-i}) \in G(k,V).$$ Ce $H$ convient.\\

$A$ est donc inversible.\\
On pose $C=A^{-1}M$ où M est le vecteur colonne 
$\begin{pmatrix}
0\\
\vdots \\
0\\
1\\
\end{pmatrix}.$\\

On pose enfin pour $j \in \{ 0,\dots,r \}$, $\epsilon_j=C_{j+1,1}$.\\
On a maintenant:\\
$$ \forall P \in G(r,V) \quad \check{\hat{f}}_{P_0}(P)=\displaystyle{\sum_{r=0}^{k}}b_{dim(P\cap P_0),r}\epsilon_r \quad ,$$\\
et donc 
$$\check{\hat{f}}_{P_0}(P)=
\left\{
\begin{array}{cc}
1 & \text{si } P=P_0,\\
0 & \text{sinon}.\\
\end{array}
\right.
$$
i.e. $\check{\hat{f}}_{P_0}=f_{P_0}.$
Le morphisme $\phi$ est donc bien injectif.\\

L'assertion $(ii)$ se démontre exactement de la même manière.$\Box$ 

\subsection{Variétés de drapeaux sur un corps fini} \label{fini2}
~\\
On conserve les notations de la section précédente.\\

On commence par rappeler la proposition suivante :

\begin{pro} \label{lemmegroupe2}
~\\
Soit $G$ un groupe, $E$ et $F$ deux $G$-ensembles finis et $\pi : E \rightarrow F$ une application $G$-équivariante. On note $r$ est le nombre d'orbites de $G$ sur $F$ et $y_1,\dots ,y_r$  des représentants de chacune de ces orbites. Alors :
$$ N(G,E)=\sum_{i=1}^r N(Stab_G(y_i),\pi^{-1}(\{y_i\})).$$
\end{pro}

On peut maintenant énoncer les deux résultats du paragraphe :

\begin{theo} \label{drapeauxfinis}
~\\
Soit $(a_1,\dots ,a_r)$ une composition de $n$. On a l'égalité suivante :
$$ N(G,\Fl_{(a_1,\dots ,a_r)}(V))= N(G,\Fl_{(a_{\sigma(1)},\dots ,a_{\sigma(r)})}(V)) \text{ pour tout $\sigma$ dans $\mathfrak{S}_r$.}$$
\end{theo}

\begin{theo} \label{drapeauxfinis2}
~\\
Soit $\underline{a}=(a_1,\dots ,a_r)$ et $\underline{b}=(b_1,\dots ,b_r)$ deux compositions de $n$ telles que $P(\underline{a}) \geq P(\underline{b})$, on a : 
$$ N(G,\Fl_{\underline{a}}(V)) \leq N(G,\Fl_{\underline{b}}(V)).$$
\end{theo}

{\sl Démonstration du théorème \ref{drapeauxfinis}:}\\
Il suffit de démontrer le résultat pour $\sigma=(i,i+1)$ avec $i$ dans $\{ 1,\dots,r-1 \}$. Dans ce cas on applique la proposition \ref{lemmegroupe2} aux flèches naturelles 
$$\Fl_{(a_1,\dots ,a_r)}(V) \rightarrow \Fl_{(a_1,\dots,(a_i+a_{i+1}),\dots,a_r)}(V)$$
  et 
$$\Fl_{(a_{\sigma(1)},\dots ,a_{\sigma(r)})}(V) \rightarrow \Fl_{(a_1,\dots,(a_i+a_{i+1}),\dots,a_r)}(V)$$  pour calculer les nombres $ N(G,\Fl_{(a_1,\dots ,a_r)}(V))$ et $N(G,\Fl_{(a_{\sigma(1)},\dots ,a_{\sigma(r)})}(V))$, le résultat découle alors de la proposition \ref{grassfin}.$\Box$\\

On démontre un résultat analogue à la proposition \ref{bleu} en utilisant la même technique, et on en déduit toujours d'après le lemme \ref{comb} le théorème \ref{drapeauxfinis2}.

\section{Annexe : un lemme combinatoire}
Cette section purement technique est une variation sur l'idée des "raising operators" du livre de Macdonald.\\
 
Soit $n$ un entier positif. Soit $(\lambda_1,\dots ,\lambda_r)$ dans $\mathcal{P}(n)$, et soit $i$ dans $\{ 1,\dots,r-1 \}$, on pose $R_{i}(\lambda)=P(\lambda_1,\dots ,\lambda_i+1,\lambda_{i+1}-1,\dots ,\lambda_r)$.
On a le résultat suivant :

\begin{lem} \label{comb}
~\\
Soit $\lambda=(\lambda_1,\dots ,\lambda_r)$ et $\mu=(\mu_1,\dots ,\mu_r)$ dans $\mathcal{P}(n)$. On a  $\lambda \geq \mu$ si et seulement s'il existe des entiers positifs $a_1,\dots ,a_{r-1}$ tels que $\lambda= R_{r-1}^{a_{r-1}} \circ \dots \circ R_{1}^{a_1} (\mu)$.
\end{lem}


\begin{thebibliography}{99}
\renewcommand{\baselinestretch}{0.9}

\bibitem{KAC} Kac, V.G. {\em Infinite root systems, representations of graphs and invariant theory ,II}. Journal of Algebra {\bf 78} (1982) 141-162.

\bibitem{McD} Macdonald, I.G. {\em Symmetric functions and Hall polynomials}. Oxford Mathematical monographs, 1979.

\bibitem{PYA} Pyasetskii, V. {\em Linear Lie groups acting with finitely many orbits}. Functional Anal. Appl. {\bf 9} (1975) 351-353.

\bibitem{ROS} Rosenlicht, M. {\em A remark on quotient spaces}. An.Acad.Brasil.Ci. {\bf 35}. (1963) 487-489.

\end{thebibliography}
\end{document}